
\input mssymb.tex
\magnification=\magstep1
\font\titlefont=cmbx10 scaled \magstep2
\def\F{\Bbb F}

\def\Z{\Bbb Z}

\def\P{\Bbb P}
\def\A{\Bbb A}

\def\vs{\vskip 1.0pc}
\def\ni{\noindent}
\def\ld{\ldots}

\noindent
\centerline{\titlefont On the Existence of Supersingular Curves}
\centerline{\titlefont Of  Given Genus }
\noindent
\vskip 2pc
\centerline{\titlefont Gerard van der Geer and Marcel van der Vlugt}
\vskip 2.0pc\ni
{\bf Introduction}
\vs\ni
In this note we shall show that there exist supersingular curves for every
positive genus in characteristic 2. Recall that an irreducible smooth
algebraic curve $C$ over an algebraically closed field $\F$ of characteristic
$p>0$ is called {\sl supersingular} if its jacobian is isogenous to a
product of supersingular elliptic curves. An elliptic curve is called
supersingular if it does not have points of order $p$ over $\F$. It is
not clear a priori that there exist such curves for every genus.
Indeed, note that  in the moduli space $A_g\otimes \F_p$ of principally
polarized abelian varieties the locus of supersingular abelian varieties
has dimension $[g^2/4]$ (cf. [O, L-O]), while the locus of jacobians has
dimension $3g-3$ for $g>1$. Therefore, as far as dimensions are
concerned  there is no  reason why these loci should
intersect for $g\geq 9$.

In this paper we construct for every integer $g>0$ a supersingular
curve of genus $g$ over the field $\F_2$. In particular this shows
that for every $g>0$ there exists an irreducible curve of genus $g$ whose
jacobian is isogenous to a product of elliptic curves.   We refer to
[E-S] for related questions in characteristic $0$. We do our
construction  by taking a suitable fibre product of Artin-Schreier
curves. This construction is inspired by coding theory, where the
introduction of generalized Hamming weights led us  to consider such
products, cf. [G-V 2].

More generally, we are able to construct   in
characteristic $p$  a supersingular curve over $\F_p$   of any genus $g$
whose $p$-adic expansion consists of the digits $0$ and $(p-1)/2$ only.
We can also count on how many moduli the construction depends.
\vs\ni
{\bf \S1  Fibre products of Artin-Schreier curves.} \vs\ni Let $\F$ be a
fixed algebraic closure of the prime field $\F_2$. We consider a finite
dimensional $\F_2$-linear subspace ${\cal L}$ of the function field
$\F(x)$. Define the operator $\wp$ on $\F(x)$ by $\wp (f) = f^2+f$. We
shall assume that ${\cal L}\cap \wp (\F(x))= \{ 0 \}$.

To an element $f \in {\cal L} - \{ 0 \}$ we associate the complete
non-singular Artin-Schreier curve $C_f$ with affine equation

$$y^2+y=f.$$
Choose a basis $f_1, \ld, f_k$ of ${\cal L}$ and let $\phi_i: C_{f_i} \to \P^1$
be the morphism given by the inclusion $\F(x) \subset \F(x,y)$. Then we define
a curve $C^{\cal L}$ by
$$C^{\cal L}= { \hbox {  \rm Normalization of }
}(C_{f_1}\times \ld \times C_{f_k}),$$
where the product means   the fibre product taken with respect
to the morphisms $\phi_i$. Up to $\F(x)$-isomor\-phism the curve $C^{\cal
L}$ is independent of the chosen basis of ${\cal L}$.

In the following we need some properties of the curve  $C^{\cal L}$; the
reader can find a proof in [G-V 2].
\vs\ni
\proclaim (1.1) Proposition. i) The jacobian of $C^{\cal L}$
decomposes up to isogeny as
$${\rm Jac}(C^{\cal L})\thicksim \prod_{f \in {\cal L}- \{ 0 \} }{\rm
Jac}(C_f)$$ and therefore the genus $g(C^{\cal L})$ can be expressed as
$$g(C^{\cal L})= \sum_{f \in {\cal L} - \{ 0 \} }g(C_f)$$
in terms of the genera of the $C_f$. \par
\smallskip
\proclaim (1.2) Corollary. Suppose that for all $f \in {\cal L}- \{ 0
\}$ the curve $C_f$ is supersingular or rational. Then the fibre
product $C^{\cal L}$ is supersingular or rational. \par \vs\ni As
ingredients for our fibre product we shall use special Artin-Schreier
curves. We consider for $h \geq 1$ the vector space ${\cal R}_h$ of
$2$-linearized polynomials $$\{ R= \sum_{i=0}^h a_ix^{2^i} : a_i \in \F
\}$$ and define
$${\cal R}_h^* = \{ R \in {\cal R}_h : a_h \neq 0 \}.$$ We proved in
[G-V 1] the following result.
\smallskip\ni
\proclaim (1.3) Proposition. The Artin-Schreier curve $C_R$ with affine
equation $y^2+y=xR(x)$ for $R \in {\cal R}_h^*$ is a
(hyperelliptic) supersingular curve of genus $2^{h-1}$. \par
\vs\ni
{\bf \S 2 The Construction.}
\vs\ni
In this section we describe how to construct a curve of a given genus
in characteristic 2. Here the construction is done over a
finite extension of the prime field. In Section 3 we shall
show that we can find such a curve over the prime field $\F_2$.
\smallskip\ni
\proclaim (2.1) Theorem. Let $\F$ be an algebraic closure of $\F_2$.
For every positive genus $g$ there exists a supersingular curve over
$\F$. \par \smallskip\ni {\bf Proof.} Take $g>0$ and write $g$ as a
dyadic expansion in the form $$g=2^{s_1}(1+ \ld +2^{r_1})+2^{s_2}(1+\ld +
2^{r_2})+ \ld + 2^{s_t}(1+ \ld + 2^{r_t}),\eqno(1)$$ where $s_i, r_i \in
\Z_{\geq 0}$ and $s_i\geq s_{i-1}+r_{i-1}+2$ for $i=2, \ld, t$. We now
choose for $i=1, \ld, t$ a $\F_2$-linear subspace ${L}_i$ of $\F(x)$
contained in ${\cal R}_{u_i}^*\cup \{ 0\} $ with $u_i=(s_{i}+1)-
\sum_{j=1}^{i-1}(r_j+1)$ and $\dim ({L}_i)=r_i+1$. We put ${\cal
L}=\oplus_{i=1}^t (x{L}_i)$. It follows directly from Propositions 1 and
2 that $C^{\cal L}$ is supersingular and that since $u_{i+1} \geq u_i+1$
for $1\leq i \leq t-1$ the genus satisfies $$\eqalign{ g(C^{\cal L})=&
\sum_{f \in {\cal L} - \{ 0 \} } g(C_f)\cr =&\sum_{i=1}^t2^{u_i-1}\cdot
2^{\sum_{j=1}^{i-1} (r_j+1)}(2^{r_i+1}-1) \cr =&\sum_{i=1}^t
2^{s_i}(2^{r_i+1}-1).\cr}$$ This last expression yields the expression
for $g$ in (1), hence $g(C^{\cal L})=g$. $\square$ \vs\ni
{}From the preceding proof we conclude that there exists supersingular
curves of genus $g>0$ already over the field $\F_{2^m}$ with
$m=\max_{1 \leq i \leq t}(r_i+1)$, where the $r_i$
occur in the expansion (1).
\vs\ni
{\bf Example.} We construct a supersingular curve of genus 30. We write
$30=2(1+2+2^2+2^3)$ and this tells us that $t=1$, $s_1=1$ and $r_1=3$.
So our curve is defined over the finite field
$\F_{16}$. We set $\F_{16}= \F_2(\alpha)$ with $\alpha^4+\alpha+1=0$.
Let $L\subset {\cal R}_{u_1}^*\cup \{ 0 \}= {\cal R}_2^*\cup \{ 0 \}$ be
the $4$-dimensional space generated by $x^4, \alpha x^4, \alpha^2 x^4$
and $\alpha^3x^4$.  Then ${\cal L}= xL$ and $C^{\cal L}$ is the desired
supersingular curve of genus 30. Its function field
${\cal F}=\F_{16}(x,y_0,y_1,y_2,y_3)$  with $y_i^2+y_i=\alpha^ix^5$ is a
Galois extension of $\F_{16}(x)$ of degree $16$. Consider the element
$y=\sum_{i=0}^3 \alpha^iy_i$. Then for all non-trivial $\sigma \in {\rm
Gal}({\cal F}/ \F_{16}(x))$ we have $\sigma(y)\neq y$, hence ${\cal
F}=\F_{16}(x,y)$. We obtain  $$\eqalign{y^{16}+y&=\sum_{i=0}^3
\alpha^i(y_i^{16}+y_i)\cr &= \sum_{i=0}^3
\alpha^i(\alpha^{8i}x^{40}+\alpha^{4i}x^{20}+
\alpha^{2i}x^{10}+\alpha^ix^5)\cr
&=\alpha^{6}x^{40}+x^{20}+\alpha^{12}x^{10}+ \alpha^9x^5.\cr}$$
Thus we have found a supersingular curve of genus $30$ over the field
$\F_{16}$ with affine equation
$$y^{16}+y=\alpha^{6}x^{40}+x^{20}+\alpha^{12}x^{10}+ \alpha^9x^5.$$

\vs\ni
{\bf \S3 Equations for fibre products}
\vs\ni
The preceding example suggests to study curves of the following type. We
consider curves $C=C_{S,R}$ defined by an equation
$$S(y)=xR_1(x)+(xR_2(x))^2+\ld + (xR_{n}(x))^{2^{n-1}},\eqno (2)$$ where
$S \in \F[y]$ is a $2$-linearized polynomial
$S=y^{2^{n}}+A_{n-1}y^{2^{n-1}}+ \ld + A_0y$ with $A_0\neq 0$ and where
the $R_i \in \F[x]$ for $i=1,,\ld, n$ are also 2-linearized polynomials
(not all $0$). We shall assume for a moment that this equation defines
an irreducible curve. Consider the $\F_2$-vector space    $$\Sigma= \{
\sigma \in \F : S(\sigma) = 0 \}.$$  An element $\sigma \in \Sigma$ acts
on $C$ via $y \mapsto y+\sigma$. Thus the curve $C$ is a Galois covering
of $\P^1$ with Galois group of type $\Sigma \cong (\Z/ 2 \Z)^n$. A
$\F_2$-linear subspace $\Sigma'$ of $\Sigma$ of codimension $1$ defines
an irreducible quotient curve $C/ \Sigma'$. If $\sigma \in \Sigma -
\Sigma'$ then the linear subspace $\Sigma'$  corresponds to a splitting
of  $$ S= B(B+B(\sigma )),\eqno(3)$$ where  $B$ is the 2-linearized
monic polynomial of degree $2^{n-1}$ in $\F[y]$ with zero set $\Sigma'$.
Note that $B(\sigma)\in \F^*$ is independent of the choice of $\sigma
\in \Sigma - \Sigma'$. If we put$$B=y^{2^{n-1}} + B_{n-2}y^{2^{n-2}}+
\ld + B_0y, \eqno(4) $$ and $\beta=B(\sigma)$ then by comparing
coefficients, (3) is equivalent to the system of equations
$$\eqalign{\beta B_0+0&= A_0,\cr
B_{i-1}^2+\beta B_i&=A_i\quad {\rm for } \quad i=1,\ld,n-2,\cr
B_{n-2}^2+\beta &= A_{n-1}.\cr}\eqno(5) $$
The compatibility of (5) comes down to the equation
$$\sum_{j=1}^n {A_{n-j}^{2^{j-1}}\over
\beta^{2^j-1}}=1 \quad {\rm or } \quad \beta^{2^n}+\sum_{j=1}^n
A_{n-j}^{2^{j-1}}\beta^{2^n-(2^j-1)}
 =0 .\eqno(6)$$
Observe that $\alpha=\beta^{-1}$ satisfies a {\sl linearized} equation,
namely
$$A_0^{2^{n-1}}\alpha^{2^n} +
A_1^{2^{n-2}}\alpha^{2^{n-1}}+ \ld + A_{n-1}\alpha^2 + \alpha
=0 .\eqno(7)$$
 Define the $\F_2$-vector space $$A= \{ \alpha \in \F : \alpha { \hbox
{ \rm satisfies } }   (7) \} .$$
The elements of $A-\{ 0 \}$ parametrize the hyperplanes of
$\Sigma'$ of $\Sigma$. The hyperplane corresponding to
$\alpha \in A -\{0\}$ will be denoted by $\Sigma_{\alpha}$.
Moreover, we set
$$T= xR_1+(xR_2)^2+ \ld +
(xR_{n})^{2^{n-1}}.$$

\smallskip\ni
\proclaim (3.1) Lemma. Each quotient curve $C_{\alpha}:=C/
\Sigma_{\alpha}$ with  $\alpha\in A-\{ 0 \}$  is of the form   $$w^2+w=
\alpha^2T,\eqno(8) $$ where $w=\alpha B$ with $B$ corresponding to
$\Sigma_{\alpha}$.  \par \smallskip\ni {\bf Proof.} One checks that $w$
is invariant under $y\mapsto y+\sigma$ with $\sigma \in
\Sigma_{\alpha}$. Substitution of (3) in (2) yields (8). $\square$
\smallskip\ni \proclaim
(3.2) Corollary. The curve in (8) is $\F[x]$-isomorphic to $$w^2+w=
\alpha^2xR_1+\alpha xR_2 + \alpha^{2^{-1}}xR_3+\ld +
\alpha^{2^{-(n-2)}}xR_{n}. \eqno(9)$$ and therefore it is supersingular
if not rational. \par \smallskip \ni {\bf Proof.} Consider the
$\F[x]$-isomorphism $$ w \mapsto w+\sum_{i=2}^{n} \sum_{j=0}^{i-2}
(\alpha^{2^{2-i}}xR_i)^{2^j}.$$  and apply Proposition (1.3). $\square$
\smallskip\ni
\proclaim (3.3) Proposition. If the curve $C$ defined
by (2) is irreducible then it is a fibre product which is supersingular
if not rational. Its  jacobian is up to isogeny the product of the
jacobians of the curves given in (8) with $\alpha \in A - \{ 0 \}$. \par
\smallskip\ni {\bf Proof.} Choose a $\F_2$-basis of $A$, say
$\alpha_1,\ld, \alpha_n$. The  curves $C_{\alpha}=C/ \Sigma_{\alpha}$
are quotients of $C$, hence $C$ admits a morphism $\phi: C \to
C_{\alpha_1}\times \ld \times C_{\alpha_n}$, the fibre product of the
$C_{\alpha_i}$ with respect to the (canonical)  maps $C_{\alpha_i}\to
\P^1$. Since the $\alpha_i$ are $\F_2$-independent, Galois theory and
Lemma (3.1) yield that  the function fields of the curves
$C_{\alpha_1}\times \ld \times C_{\alpha_j}$  and
$C_{\alpha_{j+1}}$ are linearly disjoint for $j=1,\ld, n-1$. So the fibre
product $C_{\alpha_1}\times \ld \times
C_{\alpha_n}$ is a covering of degree $2^n$ of $\P^1$. Since $C$ is also
a covering of $\P^1$ of degree $2^n$ it follows that $\phi$  is an
isomorphism. The Proposition now follows from  Proposition (1.1)
and from Corollary (3.2). $\square$

\smallskip\ni The condition that the curve
$C$ be irreducible is given in the following Lemma.
\smallskip\ni
\proclaim
(3.4) Lemma.  The curve defined by (2) is irreducible if and only if the
$n$-dimensional $\F_2$-vector space ${\cal L}$ of functions $\alpha^2T$
with $\alpha \in A$   satisfies ${\cal L} \cap \wp (\F(x))=\{ 0 \}$. \par
\smallskip\ni {\bf Proof.} The implication ``$\Rightarrow$" follows from
 Proposition 3.3. As to the implication ``$\Leftarrow$", we use the
theory of Artin-Schreier extensions (see [B], Ch.\ V, \S 11). According
to that theory the compositum of the function fields $\F(x, w_{\alpha})$
with $w_{\alpha}=\alpha B$ satisfying (8) has degree  $$\# ({\cal L}/
{\cal L} \cap \wp (\F(x)))=2^n$$ over $\F(x)$. Comparison  with the
degree of $y$ in (2) shows the irreducibility. $\square$
\smallskip
\noindent
\proclaim (3.5) Theorem. For every  integer $g>0$ there exists a
supersingular curve of genus $g$ over the prime field $\F_2$. \par
\smallskip
\noindent
{\bf Proof.} We construct a supersingular curve of the form (2) with
prescribed genus $g>0$. Recall the binary expansion of $g$ given in (1)
$$g=2^{s_1}(1+ \ld +2^{r_1})+2^{s_2}(1+\ld + 2^{r_2})+ \ld + 2^{s_t}(1+
\ld + 2^{r_t})$$ where $s_i, r_i \in \Z_{\geq 0}$ and $s_i\geq
s_{i-1}+r_{i-1}+2$ for $i=2, \ld, t$. By $w$ we denote the binary
weight $w=\sum_{i=1}^t (r_i+1)$ of $g$. First we determine the LHS
$S(y)\in \F_2[y]$ of (2) and the $w$-dimensional $\F_2$-vector space $A$.

We start with $F_0(x)=x$ and we construct inductively a sequence of
$\F_2$-linearized polynomials $F_i \in \F_2[x]$ for $i=1, \ld , t$ as
follows. We set
$$F_i= (F_{i-1})^{2^{r_i+1}}+F_{i-1}.$$
Obviously, $F_{i-1}$ divides $F_i$ for $i=1, \ld , t$.

Let $S(y)= F_t(y) \in \F_2[y]$. It has degree $2^w$. Furthermore we
define $\F_2$-linear spaces
$$A^{(i)}= \{ \alpha \in \F : F_i(\alpha)=0 \}\quad {\rm for} \quad
i=1,\ld ,t.$$We set $A=A^{(t)}$. By the divisibility property of the
$F_i$ the subspaces $A^{(i)}$ form a flag in $A$.

Now we consider for $1\leq i \leq t-1$ the polynomials
$$F_i(\alpha, x)= F_i(\alpha)x^{2^{h_i}+1} \in \F_2[\alpha, x],$$
where the $h_i= s_{i+1}-w+1+\sum_{j=i+1}^t (r_j+1)$ form a
monotonically increasing sequence.
We define for $j=0, \ld, w-1$ polynomials $xR_{w-j}(x) \in \F_2[x]$ with
$R_{w-j}$ 2-linearized by writing
 $$\sum_{i=1}^{t-1} F_i(\alpha, x)=
\sum_{j=0}^{w-1} xR_{w-j}(x)\alpha^{2^j}.$$
Here $xR_{w-j}$ is the sum (possibly empty) of those monomials
$x^{2^{h_i}+1}$ occuring in the polynomials $F_i(\alpha, x)$ which have
the monomial $\alpha^{2^j}$ as coefficient.

For $\alpha \in A-\{ 0 \}$
the curves $C_{\alpha}$ with equation (9) can be written as
$$w^2+w=F_{t-1}(\alpha)x^{2^{h_t}+1}+\ld + F_0(\alpha)x^{2^{h_1}+1}
\eqno(10)$$ (after we have converted the coefficients $\alpha^2, \alpha,
\ld, \alpha^{2^{-(w-2)}}$ to $\alpha^{2^{w-1}}, \alpha^{2^{w-2}},\ld,
\alpha$).  For the $2^{w-(r_t+1)}(2^{r_t+1}-1)$ values of $\alpha \in A -
A^{(t-1)}$ the irreducible Artin-Schreier curve $C_{\alpha}$ with
equation (10) has genus $2^{s_t-(w-(r_t+1))}$ and these curves
$C_{\alpha}$ contribute $2^{s_t}(1+2+\ld + 2^{r_t})$ to the genus
of the fibre product (2). The curves $C_{\alpha}$ with $\alpha \in
A^{t-1}-A^{t-2}$ contribute $2^{s_{t-1}}(1+2+\ld + 2^{r_{t-1}})$ to the
genus. Continuing in this way we see that the supersingular curve over
$\F_2$ with affine equation $$S(y)=xR_1+(xR_2)^2 + \ld +
(xR_w)^{2^{w-1}}$$ has the prescribed genus. $\square$
\smallskip
\noindent
{\bf Example.} Take $g=221=1+2^2(1+2+4)+2^6(1+2)$. We have $s_1=0,
s_2=2, s_3=6; r_1=0, r_2=2, r_3=1$ and $w=6$. We find
$$\eqalign{F_0(x)&=x,\quad
F_1(x)=x^2+x,\quad F_2(x)=x^{16}+x^8+x^2+x,\cr
F_3(x)&=x^{64}+x^{32}+x^{16}+x^4+x^2+x.\cr}$$The space $A$ equals $\{
\alpha \in \F : \F_3(\alpha)=0\}$. For $i=0, \ld , 2 $ the polynomials
$F_i(\alpha, x)$ are $$F_0(\alpha, x)=\alpha x^3,\quad
F_1(\alpha ,x)=(\alpha^2+\alpha)x^5,\quad
F_2(\alpha, x)=(\alpha^{16}+\alpha^8+\alpha^2+\alpha)x^9.$$
{}From  the identity $$\sum_{i=0}^2 F_i(\alpha,
x)=\sum_{j=0}^{w-1}xR_{w-j}(x)\alpha^{2^j}$$
we get
$$xR_6=x^9+x^5+x^3,\quad
xR_5=x^9+x^5,\quad
xR_3=x^9,\quad
xR_2=x^9,\quad xR_1=xR_4=0.$$
This gives a supersingular curve of genus $221$ defined by
$ F_3(y)=\sum_{k=1}^6
(xR_k)^{2^{k-1}}$ i.e. by the equation
$$y^{64}+y^{32}+y^{16}+y^{4}+y^2+y=
 x^{288}+x^{160}+x^{144}+x^{96}+x^{80}+x^{36}+x^{18}.$$

 \vs\ni
{\bf \S4 Number of moduli}
\vskip 1.0pc
\ni
 Here we count the number of moduli for our families. In the
investigation of the curves $y^2+y=xR(x)$ for $R \in {\cal R}_h^*$ in
[G-V 1] the polynomial $$E_{h,R}(x)=R(x)^{2^h}+ \sum_{i=0}^h
(a_ix)^{2^{h-i}}$$ of degree $2^{2h}$  played an important role.
We define the {\sl radical} ${\overline W}_R$  of $R$ as the
subspace of $\F$ formed by the elements  satisfying the equation
$E_{h,R}(x)=0$.

\smallskip\ni \proclaim (4.1) Proposition. Let
$h\geq 2$ and let $R=\sum_{i=0}^h a_ix^{2^i}$ and $R'=\sum_{i=0}^h a_i'
x^{2^i}$ be elements of ${\cal R}_h^*$. Then the curves $C_R$ and
$C_{R'}$ are isomorphic over $\F$ if and only if there exists a $\rho
\in \F^*$ such that $a_i'=a_i\rho^{2^i+1}$ for $i=1,\ld, h$. \par
\smallskip\ni
{\bf Proof.} Since both $C_R$ and $C_{R'}$ are hyperelliptic curves an
isomorphism $\alpha $ induces an isomorphism $\alpha' : \P^1 \to \P^1$
which fixes the (unique) branch point $\infty$ and is of the form $x
\mapsto \lambda x + \mu$ with $\lambda, \mu \in \F, \lambda \neq 0$. Let
${\overline W}_R$ (resp.\ ${\overline W}_{R'}$) be the radical of $R$
(resp.\ of $R'$).  Then by [G-V 1] we have $\lambda^{-1} {\overline
W}_R={\overline W}_{R'}$. This implies that $E_{h,R}(\lambda
X)=c_{\lambda}E_{h,R'}(X)$. By writing $X^{2^h}E_{h,R}(X)=\sum_i
(U_i+U_i^{2^i})$ with $U_i=a_i^{2^{h-i}}X^{2^h+2^{h-i}}$ we  see that
$$\lambda^{2^h}c_{\lambda}a_i'^{2^{h-i}}=\lambda^{2^h+2^{h-i}}a_i^{2^{h-i}}
\quad { \hbox {\rm for } } i\geq 1 \eqno(11)$$
and
$$ \lambda^{2^h}c_{\lambda}\in \F_{2^i}^*\quad {\hbox{\rm for all } }
i \geq 1 \quad{\hbox{\rm with}}\quad a_i\neq 0.\eqno(12)$$
Relation (12) implies $\lambda^{2^h}c_{\lambda} \in \F_{2^d}^* $ with
$d={\rm g.c.d.}\{ i \geq 1 : a_i \neq 0 \}$. There exists an element
$\eta \in \F_{2^d}^*$ such that we can write
$$\lambda^{2^h}c_{\lambda}= \eta^{(2^i+1)2^{h-i}} \quad {\rm for} \quad
i\geq 1\quad {\rm with }\quad a_i\neq 0.\eqno(13)$$
Substituting (13) in (11) we find
$$a_i'=(\lambda/\eta)^{2^i+1}a_i\quad {\rm for } \quad i=1,\ld , h.$$
Conversely, the relation $a_i'=\rho^{2^i+1}a_i$ for $i=1,\ld,h$ shows
that $\rho x R(\rho x)=xR'(x)+(\rho^2a_0+a_0')x^2$. Since for fixed $R
\in {\cal R}^*$ and varying $a \in \F$ the curves $C_{R+ax^2}$ are
mutually isomorphic over $\F$ we conclude $C_R \cong C_{R'}$. $\square$
\smallskip
\noindent
{\bf Remark.} The conclusion of the Lemma still holds for
$h=1$ if we restrict to isomorphisms which are isomorphisms of
Artin-Schreier coverings of $\P^1$ of degree $2$.

\smallskip\ni
\proclaim (4.2) Corollary. Let $n\geq 1$. The intersection of
the supersingular locus  with the hyperelliptic locus in the moduli
space ${\cal M}_{2^n}\otimes \F_2$ of curves of genus $g=2^n$  has
dimension $\geq n$. \par
\smallskip
Furthermore,  consider two $n$-dimensional $\F_2$-subspaces $L$ and $L'$
of polynomials $R=\sum_{i=1}^h a_ix^{2^i} \in \F[x]$ with $a_1 \neq 0$ if
$R\neq 0$. Let ${\cal L}= xL$ and ${\cal L}'=xL'$ and set $C=C^{\cal L}$
and $C'=C^{{\cal L}'}$. We then have:   \smallskip\ni
\proclaim (4.3) Lemma. The curves $C$ and $C'$ are isomorphic as
Galois covers of type $(\Z / 2\Z)^n$ of $\P^1$  if and only if there
exists a $\rho \in \F^*$   such that under $x
\mapsto \rho x$  the space ${\cal L}$ is transformed into
${\cal L}'$. \par
\smallskip\ni
{\bf Proof.} If the curves $C$ and $C'$ are
isomorphic as Galois covers of $\P^1$ then the corresponding
quotient curves $C_R$ and $C_{R'}$ of genus $>0$ (for $R\in
 L$) are isomorphic as covers of $\P^1$. By  Lemma (4.1) and the
subsequent remark this happens only if there is a $\rho \in \F^*$ such
that the transformation $x \mapsto \rho x$  transforms $xR$ into
$xR'$ for $R \in L$.  $\square$

\vs\ni
\proclaim (4.4) Proposition. Let $g\geq 2$ be written as in (1). Then the
supersingular locus in ${\cal M}_g\otimes \F_2$ has dimension $\geq
\sum_{i=1}^t (r_i+1)u_i  -1$, where $u_i=
(s_i+1)-\sum_{j=1}^{i-1}(r_j+1)$. \par \smallskip \ni {\bf Proof.} We
consider curves of the form $C^{\cal L}$ of genus $g \geq 2$ with
${\cal L}= \oplus xL_i$ as in the proof of Theorem (2.1). Let $m=
\sum_{i=1}^t (r_i+1)u_i$. Then the $m$ coefficients of the polynomials
in a basis of ${\cal L}$ which is a union of the bases of the
$(r_i+1)$-dimensional summands $xL_i$ in ${\cal R}_{u_i}^*$ define an
open subset $Q$ of affine $m$-space $\A_{\F}^m$. Take the
$(m-1)$-dimensional quotient of $Q$ under the action of $\F^*$. A given
curve $C$ of genus $>1$ can be written only in finitely many ways as a
Galois cover of $\P^1$ since $\# {\rm Aut}(C)(\F) <\infty$. Then Lemma
(4.3)  implies  that the natural morphism $Q
 \to {\cal M}_g\otimes \F_2$
is quasi-finite (onto its image). This proves our result. $\square$
\smallskip \ni
 \vs\ni  {\bf  References} \vs\ni [B]
Bourbaki : Alg\`ebre, Chapitres 4 \`a 7. Masson, Paris 1981.
\smallskip
\noindent
[E-S] T.\ Ekedahl, J-P.\ Serre: Exemples de courbes alg\'ebriques \`a
jacobienne compl\`ete\-ment d\'ecomposable. {\sl C.R.\ Acad.\ Sci.\
Paris \bf 317} (1993), 509-513.
 \smallskip
\noindent
[G-V 1] G.\ van der Geer, M.\ van der Vlugt: Reed-Muller codes
and supersingular curves I. {\sl Compositio Math. \bf 84} (1992),
333-367.
\smallskip \noindent
[G-V 2]  G.\ van der Geer, M.\ van der Vlugt:
Fibre products of Artin-Schreier curves and generalized Hamming weights
of codes. Report 93-05, University of Amsterdam (1993).
\smallskip
\noindent
[L-O] K-Z.\ Li, F.\ Oort: Moduli of  supersingular abelian
varieties. Preprint 824. Universiteit Utrecht, (1993).
\smallskip
\noindent
[O] F.\ Oort: Moduli of abelian varieties and Newton polygons. {\sl
C.R.\ Acad.\ Sci.\ Paris \bf 312} (1991), 385-389.
 \vskip 1.0pc
\noindent February 15, 1994. \vskip 1.0pc  \line {Gerard van der Geer
\hfil Marcel van der Vlugt} \line {Faculteit Wiskunde en Informatica
\hfil Mathematisch Instituut} \line{Universiteit van Amsterdam \hfil
Rijksuniversiteit te Leiden} \line{Plantage Muidergracht 24 \hfil Niels
Bohrweg 1} \line{1018 TV Amsterdam \hfil 2300 RA Leiden}
\line{The Netherlands \hfil The Netherlands}
\bye